\documentclass[aap,seceqn,MSNbibl,citesort,dvips]{arximspdf}
\usepackage{graphicx}
\doi{10.1214/10-AAP704}
\volume{21}
\issue{2}
\pubyear{2011}
\firstpage{589}
\lastpage{608}

\makeatletter

\newtheorem{theorem}{Theorem}
\newtheorem{proposition}[theorem]{Proposition}
\newtheorem{lemma}[theorem]{Lemma}

\newproclaim{assumption}{Assumption}

\newcommand{\expec}{\mathbf{E}}
\newcommand{\Oh}{{O}}
\newcommand{\im}{{i}}
\newcommand{\e}{{e}}
\newcommand{\dd}{{d}}
\newcommand{\prob}{\mathbf{P}}

\makeatother

\begin{document}
\begin{frontmatter}

\title{The Longstaff--Schwartz algorithm for L{\'e}vy models: Results
on fast
and slow convergence}
\runtitle{Longstaff--Schwartz algorithm for L{\'e}vy models}

\begin{aug}
\author[A]{\fnms{Stefan} \snm{Gerhold}\corref{}\thanksref{t2}\ead
[label=e1]{sgerhold@fam.tuwien.ac.at}\ead
[label=u1,url]{http://www.fam.tuwien.ac.at/\textasciitilde sgerhold/}}
\runauthor{S. Gerhold}
\affiliation{Vienna University of Technology}
\address[A]{Financial and Actuarial Mathematics\\
Vienna University of Technology\\
Wiedner Hauptstra\ss{}e 8--10\\
A-1040, Vienna\\
Austria\\
\printead{e1}\\
\printead{u1}} 
\end{aug}

\thankstext{t2}{Supported by the Christian Doppler Laboratory for
Portfolio Risk Management
(PRisMa Lab) \protect\url{http://www.prismalab.at/}. The Bank Austria
Creditanstalt (BA-CA) and the Austrian Federal Financing Agency
(\"{O}BFA) through CDG is gratefully acknowledged.}

\received{\smonth{3} \syear{2008}}
\revised{\smonth{3} \syear{2010}}

%
\begin{abstract}
We investigate the Longstaff--Schwartz algorithm for American option
pricing assuming that both the number of regressors and the number of
Monte Carlo paths tend to infinity. Our main results concern
extensions, respectively, applications of results by Glasserman and Yu
[\textit{Ann. Appl. Probab.} \textbf{14} (2004) 2090--2119] and
Stentoft [\textit{Manag. Sci.} \textbf{50} (2004) 1193--1203] to
several L\'evy models, in particular the geometric Meixner model. A
convenient setting to analyze this convergence problem is provided by
the L\'evy--Sheffer systems introduced by Schoutens and Teugels.
\end{abstract}

%

\pdfkeywords{62P05, 33C45, Option pricing, dynamic programming,
Monte Carlo, regression, orthogonal polynomials, Levy--Meixner systems.}

\end{frontmatter}

\section{Introduction}

PDE or tree methods for pricing financial products become ineffective
in the presence of
many stochastic factors and path dependent payoff structures. When
resorting to Monte Carlo,
early exercise features like callability or flip options pose difficulties.
Typical examples are the pricing of callable \mbox{LIBOR} exotics with the
LIBOR market model \cite{BrMe06}
or the valuation of life insurance contracts with early exercise
features \cite{BaBiMi09}.

The least squares Monte Carlo approach by Longstaff and Schwartz \cite
{LoSc01} has become the standard
method to deal with such American/Bermudan products. It proceeds by
backward induction
and estimates value functions by regression on a prescribed set of
basis functions.
The computed exercise strategy is suboptimal, resulting in a lower
bound for the option price;
see Belomestny, Bender and Schoenmakers \cite{BeBeSc09} for recent work on upper bounds.
Fouque and Han \cite{FoHa09} discuss numerical aspects of American
option pricing,
including variance reduction.

The convergence analysis of the Longstaff--Schwartz algorithm was commenced
in the original paper \cite{LoSc01} and was carried out in detail
by Cl\'{e}ment, Lamberton and Protter \cite{ClLaPr02}. They show
convergence of the regression approximation
to the true Bermudan price and convergence of the Monte Carlo procedure
for a fixed number of basis functions.
Glasserman and Yu \cite{GlYu04} and Stentoft \cite{St04} have analyzed
settings in which
the number of basis functions and the number of simulation paths
increase together. In particular, Glasserman and Yu \cite{GlYu04} have
shown that the number of paths must grow exponentially in
the number of basis functions if the underlying process is Brownian
motion or geometric Brownian motion.
On the other hand, Stentoft \cite{St04} appealed to results on series
estimators \cite{deJo02,Ne97}
to obtain polynomial growth for rather general models, assuming
that the underlying has a bounded state space.
The latter assumption was also imposed by Eglof, Kohler and
Todorovic \cite{Eg05,EgKoTo07} in the analysis of
their extension of the Longstaff--Schwartz algorithm.

In the present paper we discuss the applicability of Stentoft's results
to exponential L\'evy models
and extend Glasserman and Yu's analysis to several models, including
the Meixner model \cite{Gr99,Sc02}.
These latter results provide an application of the neat martingale
properties that Schoutens
and Teugels \cite{Sc00,ScTe98} found for certain
L\'evy processes and families of orthogonal polynomials.

In the following section we recall the dynamic programming principle
and the Longstaff--Schwartz algorithm.
We show how Stentoft's \cite{St04} convergence result can be applied to
L\'evy models, in particular,
to the Meixner model. This involves discussing the assumption of a
bounded underlying and the smoothness
of the value functions occurring in the backward induction.

In Section \ref{se:gy} we describe the problem that Glasserman and
Yu \cite{GlYu04} treated. The
main difference to Stentoft's setting is the unbounded support of the
underlying.
Section \ref{se:meixner} recalls the notions of Sheffer system and L\'
{e}vy--Meixner system.
Besides Brownian motion, this theory yields four processes that lend
themselves to the investigation:
the Meixner, standard Poisson, Gamma and Pascal processes \mbox{\cite
{DuGeSh91,GeSh94}}.
In Section \ref{se:single} we assume that our option has only three
exercise opportunities resulting
in a single regression and show
how fast the number of simulation paths must increase in order to
ensure convergence of the Longstaff--Schwartz algorithm
for a growing number of basis functions.
Finally, Section \ref{se:multi} contains an analogous bound for the
multi-period setting, which is weaker,
but upon inversion still leads to the same critical asymptotic rate as
the single-period case.
In the course of the proofs it turns out that the different critical
rate pertaining to Brownian motion
stems from the comparatively slow growth of the linearization
coefficients of the associated
L\'evy--Meixner system, namely, the Hermite polynomials.

\section{Bounded state space and fast convergence}\label{se:stentoft}

Suppose that our asset follows a Markov process $S_t$.
We assume throughout the paper that the interest rate is zero; extending
our results to a constant interest rate $r>0$ is trivial.
Consider a Bermudan option (which may serve as a proxy for an American
option) that can be exercised at the times
$0= t_0<\cdots< t_m$. The payoff from exercise is $h_n(S_{t_n})$ for
given functions $h_n$, $0\leq n\leq m$.
By the dynamic programming principle the option value at time $t_0=0$
equals $V_0 = \max\{h_0(S_{0}), C_0(S_{0})\}$,
where the continuation values $C_n$ are given by
\begin{eqnarray*}
C_m(x) &=& 0, \\
C_n(x) &=& \expec[\max\{h_{n+1}(S_{t_{n+1}}), C_{n+1}(S_{t_{n+1}})\}
\mid S_{t_n} = x],\qquad  0\leq n< m.
\end{eqnarray*}
Suppose that $N$ sample paths of the underlying are simulated.
Longstaff and Schwartz \cite{LoSc01} propose to approximate the
continuation values by a linear combination
of basis functions $\psi_{nk}$,
\[
C_n(x) \approx\sum_{k=0}^K \beta_{nk} \psi_{nk}(x) = \beta
_n^{\mathrm
{T}} \psi_n(x),
\]
where $\beta_n = (\beta_{n0},\ldots, \beta_{nK})^{\mathrm{T}}$ is a
vector of real numbers which is estimated
by regression over the simulated paths and
$\psi_n(x) = [\psi_{n0}(x),\ldots,\psi_{nK}(x)]^{\mathrm{T}}$.

To obtain a good convergence result as $N$ and $K$ both tend to
infinity, Stentoft \cite{St04} assumes that
samples above and below certain thresholds are discarded.
So let us fix finite truncation intervals $I_1,\ldots,I_m \subset
{]0,\infty[}$ and discard
all sample paths with $S_{t_n} \notin I_n$ when estimating the
continuation value $C_n(x)$.
We are then estimating the following ``truncated'' continuation values:
\begin{eqnarray*}
C_m^{\mathrm{tr}}(x) &=& 0, \\
C_n^{\mathrm{tr}}(x) &=& \mathbf{1}_{I_n}(x) \cdot\expec[\max\{
h_{n+1}(S_{t_{n+1}}), C^{\mathrm{tr}}_{n+1}(S_{t_{n+1}})\} \mid S_{t_n}
= x], \qquad
1\leq n< m, \\
C_0^{\mathrm{tr}}(x) &=& \expec[\max\{h_{1}(S_{t_{1}}), C^{\mathrm
{tr}}_{1}(S_{t_{1}})\} \mid S_{0} = x].
\end{eqnarray*}
The option value at time $t_0=0$ approximately equals
%
\begin{equation}\label{eq:V0tr}
V_0^{\mathrm{tr}} = \max\{h_0(S_{0}), C_0^{\mathrm{tr}}(S_{0})\}.
\end{equation}
Outside of the truncation intervals $I_n \subset{]0,\infty[}$ we
extrapolate by zero
since it does not matter in the theoretical analysis.

Besides truncation, another possibility to make the state space bounded
would be absorption of the underlying process
at some lower and upper bounds \cite{Eg05,EgKoTo07}. This, however,
causes atoms in the distribution
so that Stentoft's result is no longer applicable as it requires the
existence of a density.

We assume in the present section that the underlying has the following dynamics.
(Recall that we suppose throughout that the interest rate is zero.)
\renewcommand{\theassumption}{\Alph{assumption}}
\begin{assumption}[(Exponential L\'evy dynamics)]\label{assA}
The risk neutral dynamics of the underlying are
\[
S_t = S_0 \exp( X_t),
\]
where $X_t$ is a L\'evy process with $X_0=0$.
The support of $X_t$ is the whole real line for $t>0$ and $X_t$ has a continuous
density function.
\end{assumption}
\begin{assumption}[(Value smoothness)]\label{assB}
Let the function $h$ be of, at most, linear growth and such that
$h(S_T)$ is integrable for each $T>0$. Then $\expec[h(S_T) \mid S_0 =
x]$ is a $C^1$-smooth function of $x$.
\end{assumption}

Without going into detail we note
that Stentoft \cite{St04} imposes the following additional assumptions:
\begin{assumption}[(Further technical assumptions)]\label{assC}
The basis functions are shifted Legendre polynomials,
the continuation values $C_n(S_{t_n})$ are in the $L^2$-span of the regressors,
the simulated paths are independent and the probability that the
exercise payoff exactly
equals the continuation value is zero.
\end{assumption}

Now Stentoft's main result (\cite{St04}, Theorem 2),
specialized to L\'evy models, reads as follows.
(By ``truncated algorithm'' we mean that we discard the samples outside
the intervals $I_n$
as explained above.)
\begin{theorem}\label{thm:stentoft}
Fix arbitrary finite truncation intervals $I_1,\ldots,I_m$ contained in
${]0,\infty[}$ and assume that
Assumptions \ref{assA}--\ref{assC} hold.
Let $N$ (the number of paths) and $K$ (the number of basis functions)
tend to infinity such that $K^3/N\to0$. Then the option prices
computed by the truncated Longstaff--Schwartz algorithm converge
to $V_0^{\mathrm{tr}}$,
defined by (\ref{eq:V0tr}).
\end{theorem}

If the truncation intervals are large enough, then one would hope that
the approximate price $V_0^{\mathrm{tr}}$ is close
to the exact price $V_0$. We will now show that this is indeed the case
for L\'evy models,
assuming mild integrability and (at most) linearly growing payoff functions.
\begin{assumption}[(Integrability)]\label{assD}
For each $t$ there are $p>1$ and $p'>0$ such that $S_t^p$ and
$S_t^{-p'}$ are integrable.
\end{assumption}
\begin{assumption}[(Linear payoff growth)]\label{assE}
The payoff functions grow at most linearly,
%
\begin{equation}\label{eq: bound h}
|h_n(x)| \leq c(1+x), \qquad  x\geq0, 1\leq n \leq m \mbox{, for
some } c > 0.
\end{equation}
\end{assumption}
%
\begin{theorem}\label{thhh}
Assume that Assumptions \ref{assA}, \ref{assD} and \ref{assE} are satisfied and that the
truncation intervals satisfy
\[
I_n = [b_n^{-1}, b_n], \qquad 1\leq n < m,
\]
where
%
\begin{equation}\label{eq:b expo}
b_n = b_{n+1}^\nu, \qquad 1\leq n< m-1,
\end{equation}
with
\[
\nu= \min\biggl\{\frac{p'}{p'+q}, \frac{p}{p+q}\biggr\}  \quad\mbox{and}\quad
\frac1p + \frac1q = 1.
\]
Then $V_0^{\mathrm{tr}}$ converges to the exact option price $V_0$
as $b_m$ tends to infinity.
\end{theorem}

Note that $\nu= 1-1/p$ in (\ref{eq:b expo}) if $p=p'$. In particular,
if all moments of the underlying
and its reciprocal exist, like in the Black--Scholes model, then
the exponent $\nu\in{]0,1[}$ is arbitrarily close to $1$.
\begin{pf*}{Proof of Theorem \ref{thhh}}
A trivial induction, using the martingale property of~$S_t$, shows
that the continuation values $C_n(x)$
and $C_n^{\mathrm{tr}}(x)$ satisfy the bound (\ref{eq: bound h}) too.
We will show that for all $n$
\[
C^{\mathrm{tr}}_n(x) = C_n(x) + o(1)  \qquad\mbox{as } b_m
\to\infty\mbox{, uniformly w.r.t. } x\in I_n.
\]
For $x\in I_n$, we have
%
\begin{eqnarray}
&& C_n(x) - C^{\mathrm{tr}}_n(x) \\
&&\qquad= \expec\bigl[\mathbf{1}_{\{S_{t_{n+1}} \notin I_{n+1}\}}
\max\{ h_{n+1}(S_{t_{n+1}}), C_{n+1}(S_{t_{n+1}}) \} \mid S_{t_n} = x\bigr]
\nonumber\\
&&\qquad\quad{} + \expec\bigl[\mathbf{1}_{\{S_{t_{n+1}} \in I_{n+1}\}}
\max\{ h_{n+1}(S_{t_{n+1}}), C_{n+1}(S_{t_{n+1}}) \} \mid S_{t_n} = x\bigr]
\nonumber\\
\label{eq:c-ctr}
&&\qquad\quad{} - \expec\bigl[\mathbf{1}_{\{S_{t_{n+1}} \in I_{n+1}\}}
\max\{ h_{n+1}(S_{t_{n+1}}), C^{\mathrm{tr}}_{n+1}(S_{t_{n+1}}) \}
\mid S_{t_n} = x\bigr] \\
&&\qquad\quad{} - \expec\bigl[\mathbf{1}_{\{S_{t_{n+1}} \notin I_{n+1}\}}
\max\{ h_{n+1}(S_{t_{n+1}}), C^{\mathrm{tr}}_{n+1}(S_{t_{n+1}}) \}
\mid S_{t_n} = x\bigr]. \nonumber
\end{eqnarray}
It follows readily from the induction hypothesis that the difference of
the second and the
third term is uniformly $o(1)$ on $I_n$, as $b_{m-1}\to
\infty$.
In the following, we write $c$ for various positive constants whose
precise value
is irrelevant.
Now let us estimate the first and the last expectation on the
right-hand side of (\ref{eq:c-ctr}).
Again, for $x\in I_n$ we use H\"older's inequality and Minkowski's inequality
to see that each of them is bounded by
%
\begin{eqnarray} \label{eq:est b}\qquad
&&
\expec\bigl[\mathbf{1}_{\{S_{t_{n+1}} \notin I_{n+1}\}}
c (1 + S_{t_{n+1}}) \mid S_{t_n} = x\bigr] \nonumber\\
&&\qquad \leq c \prob[S_{t_{n+1}} \notin I_{n+1} \mid S_{t_n} =
x]^{1/q} \cdot
\expec[(1 + S_{t_{n+1}})^p \mid S_{t_n} = x]^{1/p} \nonumber\\
&&\qquad = c\prob\biggl[x\frac{S_{t_{n+1}}}{S_{t_n}} \notin
I_{n+1}\biggr]^{1/q}
\cdot\expec\biggl[\biggl(1 + x\frac{S_{t_{n+1}}}{S_{t_n}}\biggr)^p \biggr]^{1/p}
\\
&&\qquad \leq c \biggl( 1 - F_n\biggl(\frac{b_{n+1}}{x}\biggr) + F_n\biggl(\frac{1}{x
b_{n+1}}\biggr) \biggr)^{1/q}
\biggl(1 + x \expec\biggl[\biggl(\frac{S_{t_{n+1}}}{S_{t_n}}\biggr)^p
\biggr]^{1/p}\biggr) \nonumber\\
&&\qquad \leq c b_n \biggl( 1 - F_n\biggl(\frac{b_{n+1}}{b_n}\biggr) + F_n\biggl(\frac
{b_n}{b_{n+1}}\biggr) \biggr)^{1/q},\nonumber
\end{eqnarray}
where $F_n$ is the distribution function of $S_{t_{n+1}}/S_{t_n}$.
Now note that
\begin{eqnarray*}
b_n^q \biggl(1-F_n\biggl(\frac{b_{n+1}}{b_n}\biggr)\biggr)
&=& b_n^{p+q} b_{n+1}^{-p} \biggl( \frac{b_{n+1}}{b_n}\biggr)^p
\biggl(1- F_n\biggl(\frac{b_{n+1}}{b_n}\biggr)\biggr) \\
&\leq&\biggl( \frac{b_{n+1}}{b_n}\biggr)^p \biggl(1 - F_n\biggl(\frac
{b_{n+1}}{b_n}\biggr)\biggr) \\
&=& o(1),\qquad
b_{m-1}\to\infty,
\end{eqnarray*}
where the last equality follows \cite{Fe71} from $S_{t_{n+1}}/S_{t_n}
\in L^p$.
Similarly, if $G_n$ denotes the distribution function of
$S_{t_n}/S_{t_{n+1}}$, we have
\begin{eqnarray*}
b_n^q F_n\biggl(\frac{b_n}{b_{n+1}}\biggr) &=& b_n^q \biggl(1 - G_n\biggl(\frac
{b_{n+1}}{b_n}\biggr)\biggr) \\
&=& b_n^{p'+q} b_{n+1}^{-p'} \biggl( \frac{b_{n+1}}{b_n}\biggr)^{p'}
\biggl(1 - G_n\biggl(\frac{b_{n+1}}{b_n}\biggr)\biggr) \\
&=& o(1).
\end{eqnarray*}
\upqed\end{pf*}

Besides the bounded state space, a crucial assumption of Stentoft's
result (Theorem \ref{thm:stentoft})
is the smoothness of the continuation value functions. In the
Black--Scholes model, and more generally in models
where the log-price $X_t$ has a diffusion component, they are
always $C^\infty$-smooth \cite{CoTaVo04}.
The variance Gamma model is an example of a pure jump process where the
value functions are not
necessarily continuously differentiable \cite{CoTaVo04}. In the
geometric Meixner model \cite{Gr99,Sc02,ScTe98},
on the other hand, the continuation
values are smooth, as we will now show.
Consequently, Theorem \ref{thm:stentoft} is applicable to the geometric
Meixner model
(if the mild Assumptions~\ref{assC} and \ref{assD} are satisfied).
\begin{proposition}\label{prop:meixner}
Suppose that Assumptions \ref{assA} and \ref{assD} hold and that the log-price $X_t$ is
a Meixner process.
Then Assumption \ref{assB} holds.
\end{proposition}
\begin{pf}
For fixed $t>0$ the log-price $X_t$ follows the Meixner distribution
$\mathrm{Meix}(\alpha,\beta,\mu t, \delta t)$, where
$\alpha>0$, $-\pi<\beta<\pi$, $\mu>0$ and $\delta\in\mathbb
{R}$. This
means that the density of $X_t$ equals
\[
f_t(x) = \frac{(2 \cos({\beta}/{2}))^{2\delta t}}{2\pi\alpha
\Gamma
(2\delta t)}
\e^{{\beta}/{\alpha} (x-\mu t)} \biggl|\Gamma\biggl(\delta t + \im\frac
{x-\mu
t}{\alpha}\biggr)\biggr|^2
\]
and the value function for the payoff $h(S_T)$ is
%
\begin{eqnarray}\label{eq:val func}
\expec[h(S_T) \mid S_t = x] &=& \int_{-\infty}^\infty h(\e^y x)
f_{T-t}(y) \,\dd y \nonumber\\[-8pt]\\[-8pt]
&=& \int_0^\infty h(z) f_{T-t}\biggl(\log\frac{z}{x}\biggr) \,\dd z/z .\nonumber
\end{eqnarray}
By the asymptotic formulas \cite{Sc02}
\[
f_t(x) \sim c_\pm|x|^{2\delta t-1} \e^{-|x| (\pi\pm\beta)/\alpha},
\qquad  x\to\pm\infty,
\]
and the integrability Assumption \ref{assD}, we must have $(\pi+\beta)/\alpha
> 1$.
We can now differentiate the value function (\ref{eq:val func}) under
the integral sign,
justified by the\vadjust{\goodbreak} following fact: for real $u$ and natural $k$ the quantity
\[
\frac{{\partial^k}/{\partial v^k}|\Gamma(u + \im v)|}{|\Gamma
(u +
\im v)|}
\]
grows only polynomially in $v$ as $v\to\pm\infty$.
To see this start from Lerch's formula \cite{Go01}
\[
|\Gamma(u+\im v)| = \frac{\Gamma(u+1)}{\sqrt{u^2+v^2}} \prod
_{n=1}^\infty
\biggl( 1+\frac{v^2}{(u+n)^2} \biggr)^{-1/2},
\]
hence, we have
\[
\frac{{\partial}/{\partial v}|\Gamma(u + \im v)|}{|\Gamma(u +
\im v)|}
= -\frac{v}{u^2+v^2} - v \sum_{n=1}^\infty\frac{1}{(u+n)^2+v^2}.
\]
It suffices to note that $1/[(u+n)^2+v^2] \leq1/(u+n)^2$ to see that
this expression
grows only polynomially in $v$. The higher derivatives can be dealt
with by a straightforward induction.
\end{pf}

\section{Unbounded state space and slow convergence}\label{se:gy}

If we drop the assumption that the state space of our underlying is bounded,
the convergence behavior of the Longstaff--Schwartz algorithm radically changes.
(As above, we suppose that both the number of paths and the number of
basis functions
tend to infinity.)
This is illustrated by results of Glasserman and Yu \cite{GlYu04} who
showed, assuming that
the underlying follows either Brownian motion or geometric Brownian motion,
that the number of Monte Carlo paths must grow exponentially in the number
of basis functions to retain convergence.
The first and last lines of Table \ref{ta:comparison} reflect this result;
the lines in between will be established below.

\begin{table}[b]
\caption{The highest possible number of basis functions for $N$
paths}\label{ta:comparison}
\begin{tabular*}{\tablewidth}{@{\extracolsep{\fill}}lcc@{}}
\hline
\textbf{Process} & \multicolumn{1}{c}{\textbf{Basis polynomials}}
& \multicolumn{1}{c@{}}{\textbf{\#Basis functions}} \\
\hline
Geometric Brownian motion & Monomials & $\sqrt{\log N}$ \\
Meixner & Meixner--Pollaczek & $\log N/ \log\log N$ \\
Standard Poisson & Charlier & $\log N/ \log\log N$ \\
Gamma & Laguerre & $\log N/ \log\log N$ \\
Pascal & Meixner & $\log N/ \log\log N$ \\
Brownian motion & Hermite & $\log N$ \\
\hline
\end{tabular*}
\end{table}

For the reader's convenience, our notation closely follows that
of \cite
{GlYu04}.
Recall that we assume that the interest rate is $r=0$ throughout the paper.

The variant of the Longstaff--Schwartz algorithm to be analyzed
proceeds as follows.
Start with the final continuation value $\hat{C}_m = 0$
and the final option value $\hat{V}_m = h_m$. For $n=m-1,\ldots,1$
generate $N$ sample paths $\{S_{t_1}^{(i)},\ldots,S_{t_{n+1}}^{(i)}\}$,
$1\leq i\leq N$, and set
\begin{eqnarray*}
\hat{\gamma}_n &=& \frac1N \sum_{i=1}^N \hat
{V}_{n+1}\bigl(S_{t_{n+1}}^{(i)}\bigr) \psi_n\bigl(S_{t_{n}}^{(i)}\bigr), \\
\hat{\beta}_n &=& \Psi_n^{-1} \hat{\gamma}_n, \\
\hat{C}_n &=& \hat{\beta}_n^{\mathrm{T}} \psi_n, \\
\hat{V}_n &=& \max\{h_n, \hat{C}_n\}.
\end{eqnarray*}
Finally, the initial continuation value is $\hat{C}_0(S_0) = N^{-1}
\sum
_{i=1}^N \hat{V}_1(S_{t_1}^{(i)})$
from which the initial option value is estimated by $\hat{V}_0 = \max
\{
h_0(S_0), \hat{C}_0(S_0)\}$.

There are two (minor) differences to the variant of the algorithm that
we analyzed in Section \ref{se:stentoft}: first, we assume now that a
fresh set of paths
is generated for each exercise date. Second, in the present section we
will use
explicit expressions for the $(K+1)\times(K+1)$ matrix
%
\begin{equation}\label{eq:Psi}
\Psi_n = \expec[\psi_n(S_{t_n}) \psi_n(S_{t_n})^{\mathrm{T}}],
\end{equation}
which has to be estimated by its sample counterpart in general.

In the single-period case $m=2$,
the question that Glasserman and Yu \cite{GlYu04} treated is as
follows. Suppose that there is an
\textit{exact} representation
%
\begin{equation}\label{eq:repr}
h_2(S_{t_2}) = \sum_{k=0}^K \beta_k \psi_{2k}(S_{t_2}),
\end{equation}
with unknown constants $\beta_k$. This assumption is not too
restrictive; an infinite series representation
of this kind has to be assumed anyway to get convergence of the
algorithm and since we are interested in $K\to\infty$,
we can suppose that (\ref{eq:repr}) is a good approximation of the
payoff at $t_2$.
Furthermore, assume that the martingale property
%
\begin{equation}\label{eq:mart}
\expec[\psi_{2k}(S_{t_2}) \mid S_{t_1}] = \psi_{1k}(S_{t_1})
\end{equation}
holds.
(In \cite{GlYu04}, additional deterministic factors in (\ref{eq:mart})
are allowed; we chose
to absorb these into the basis functions.)
How fast may $K$ tend to infinity compared to $N$ while
assuring that the mean square error of $\beta$ tends to zero? To this
end, Glasserman and Yu \cite{GlYu04} established the bounds
%
\begin{equation}\label{eq:upper bound}
\sup_{|\beta|=1} \expec[|\beta-\hat{\beta}|^2] \leq\frac{\|
\Psi
_1^{-1} \|^2}{N}
\sum_{j=0}^K \sum_{k=0}^K \expec[\psi_{2j}(S_{t_2})^2 \psi_{1k}(S_{t_1})^2]
\end{equation}
and
%
\begin{equation}\label{eq:lower bound}
\sup_{|\beta|=1} \expec[|\beta-\hat{\beta}|^2] \geq\frac{1}{N\|
\Psi
_1 \|^2}
\sum_{k=0}^K \expec[\psi_{2K}(S_{t_2})^2 \psi_{1k}(S_{t_1})^2] -
\frac{1}{N}.
\end{equation}
Here and in what follows, \mbox{$|\cdot|$} denotes the Euclidean vector norm
and \mbox{$\|\cdot\|$} denotes
the Euclidean (or Frobenius) matrix norm.
With regard to notation, Glasserman and Yu \cite{GlYu04} call the
coefficients in (\ref{eq:repr}) $a_k$
instead of $\beta_k$; our simplified assumption (\ref{eq:mart}) makes
both their $a$ and $\beta$
equal to our $\beta$. This has to be kept in mind when comparing (\ref
{eq:upper bound})
and (\ref{eq:lower bound}) to \cite{GlYu04}, formulas (22),
respectively, (23).

The proofs of the estimates (\ref{eq:upper bound}) and (\ref{eq:lower
bound}) are short; the bulk of the
work of Glasserman and Yu \cite{GlYu04} lies in the concrete examples
(Brownian motion and geometric Brownian motion)
and in the general analysis of the multi-period case on which we will
build in Section \ref{se:multi}.

The martingale property (\ref{eq:mart}) is convenient for estimating
the expectations in the
bounds (\ref{eq:upper bound}) and (\ref{eq:lower bound}).
Another useful property is orthogonality of the basis functions.
If $S_t$ is Brownian motion, then Glasserman and Yu \cite{GlYu04} have
shown that for $N$ paths the highest $K$, for which the mean square
error tends
to zero, is roughly $\log N$. Hermite polynomials are natural basis
functions in this case.
If the underlying process is geometric Brownian motion and monomials
are used as basis functions,
then $K$ may only be as high as $\sqrt{\log N}$. In the following
sections we show that the analogous rate for
the Meixner, Poisson, Gamma and Pascal processes is in between,
namely, $\log N/\log\log N$.

\section{\texorpdfstring{L\'evy--Meixner systems}{Levy--Meixner systems}}\label{se:meixner}

A source of basis functions and processes that satisfy martingale
equalities of the type (\ref{eq:mart})
are the L\'{e}vy--Meixner systems introduced by Schoutens and
Teugels \cite{Sc00,ScTe98}. Recall that
Meixner \cite{Me34} has determined all sets of orthogonal polynomials
$Q_k(x)$ that satisfy Sheffer's condition
\[
f(z) \exp(x u(z)) = \sum_{k=0}^\infty Q_k(x) \frac{z^k}{k!}
\]
for some formal power series $f$ and $u$ with $u(0)=0$, $u'(0) \neq0$
and $f(0) \neq0$.
Schoutens and Teugels \cite{ScTe98} introduce a time parameter $t$ via
\[
f(z)^t \exp(x u(z)) = \sum_{k=0}^\infty Q_k(x,t) \frac{z^k}{k!}
\]
and show how an infinitely divisible characteristic function, and thus
a L\'{e}vy process, can be defined
by $f$ and $u$ under appropriate conditions. Building on Meixner's
characterization, five
sets of orthogonal polynomials $Q_k(X_t,t)$ and associated L\'{e}vy
processes $X_t$ are determined which
satisfy martingale equalities of the type
\[
\expec[Q_k(X_t,t) \mid X_s] = Q_k(X_s,s),\qquad  0\leq s\leq t.
\]
This furnishes the connection between Sheffer (resp., L\'
{e}vy--Meixner) systems and condition (\ref{eq:mart}).
There are five L\'{e}vy--Meixner systems constructed from Hermite
polynomials, Charlier\vadjust{\goodbreak} polynomials $C_k(x,\mu)$,
Laguerre polynomials $L_k^{(\alpha)}(x)$, Meixner
polynomials $M_k(x;\mu
,q)$ and Meixner--Pollaczek
polynomials $P_k(x;\mu,\zeta)$, respectively. The resulting L\'evy
processes $X_t$ are standard Brownian motion $B_t$,
the standard Poisson process $N_t$, the Gamma process $G_t$, the Pascal
process $P_t$ and the Meixner process $H_t$, respectively.
See Schoutens and Teugels \cite{Sc00,ScTe98} for details on all these
processes and families
of orthogonal polynomials.
%

Brownian motion is not of interest to us since the corresponding last line
of Table \ref{ta:comparison} has been established by Glasserman and
Yu \cite{GlYu04}. As for the remaining four processes,
in the light of condition (\ref{eq:mart}), the martingale
relations \cite{Sc00}
%
\begin{eqnarray}\label{eq:mart gamma}
\expec[C_k(N_t,t) \mid N_s] &=& \biggl(\frac{s}{t}\biggr)^k C_k(N_s,s),
\nonumber\\
\expec\bigl[L_k^{(t-1)}(G_t) \mid G_s\bigr] &=& L_k^{(s-1)}(G_s),
\nonumber\\[-8pt]\\[-8pt]
\expec[M_k(P_t;t,q) \mid P_s] &=& \frac{(s)_k}{(t)_k} M_k(P_s;s,q),
\nonumber\\
\expec[P_k(H_t;t,\zeta) \mid H_s] &=& P_k(H_s;s,\zeta), \nonumber
\end{eqnarray}
valid for $0<s<t$, prompt us to choose the basis functions in
Table \ref{ta:levy}.
[Note that $(t)_k= t(t+1)\cdots(t+k-1)$ is the Pochhammer symbol.]
%
\begin{table}[b]
\caption{L\'{e}vy--Meixner systems}\label{ta:levy}
\begin{tabular*}{\tablewidth}{@{\extracolsep{\fill}}lccc@{}}
\hline
\textbf{Process} & \multicolumn{1}{c}{\textbf{Notation}}
& \multicolumn{1}{c}{\textbf{Basis polynomials} $\bolds{\psi_{nk}(x)}$}
& \multicolumn{1}{c@{}}{\textbf{Parameters}} \\
\hline
Meixner & $H_t$ & $\psi_{nk}^{\mathrm{M}}(x)=P_k(x;t_n,\zeta)$ &
$0<\zeta<\pi$ \\[4pt]
Standard Poisson & $N_t$ & $\psi_{nk}^{\mathrm{P}}(x)=t_n^k
C_k(x,t_n)$ & \\[2pt]
Gamma & $G_t$ & $\psi_{nk}^{\mathrm{G}}(x)=L_k^{(t_n-1)}(x)$ & \\[2pt]
Pascal & $P_t$ & $\psi_{nk}^{\mathrm{Pa}}(x)=(t_n)_k M_k(x;t_n,q)$ &
$0<q<1$ \\
\hline
\end{tabular*}
\end{table}
When specializing the bounds (\ref{eq:upper bound}) and (\ref
{eq:lower bound}) to our examples,
we will require the orthogonality properties
%
\begin{eqnarray}
\label{eq:orth charl}
\expec[C_k(N_t,t) C_l(N_t,t)] &=& t^{-k} k! \delta_{kl},  \\
\label{eq:orth lag}
\expec\bigl[L^{(t)}_k (G_t) L^{(t)}_l (G_t)\bigr] &=& \frac{\Gamma(k+t+1)}{k!}
\delta_{kl},\\
\label{eq:orth meixner}
\expec[M_k(P_t;t,q) M_l(P_t;t,q)] &=& \frac{k!}{(t)_k q^k} \delta_{kl},
\\
\label{eq:orth meix-poll}
\expec[P_k(H_t;t,\zeta) P_l(H_t;t,\zeta)] &=& \frac{\Gamma
(k+2t)}{(2\sin
\zeta)^{2t}k!} \delta_{kl},
\end{eqnarray}
as well as a way to express the squares of the basis functions as
series of basis functions. We will denote
by $d_{ki}(t_n)$ the \textit{linearization coefficients} in the expansion
%
\begin{equation}\label{eq:conn}
\psi_{nk}(x)^2 = \sum_{i=0}^{2k} d_{ki}(t_n) \psi_{ni}(x).
\end{equation}
Where distinction is necessary, the linearization coefficients
corresponding to the four families in Table \ref{ta:levy} will be written
as $d_{ki}^{\mathrm{P}}(t_n)$, $d_{ki}^{\mathrm{G}}(t_n)$,
$d_{ki}^{\mathrm{Pa}}(t_n)$
and $d_{ki}^{\mathrm{M}}(t_n)$, respectively.
The same superscripts will adorn other quantities to distinguish the
four cases, namely, the
Meixner, Poisson, Gamma and Pascal process as in Table~\ref{ta:levy}.

Among these processes, the Meixner process has the most significance in
applications. Clearly, a financial model
will impose geometric Meixner dynamics (as in Proposition \ref
{prop:meixner}) rather than
the linear process which may become negative. But then a convergence
analysis in the spirit of Glasserman and Yu \cite{GlYu04}
is impossible with polynomial basis functions as the geometric Meixner
process does not have finite moments of all orders.
Instead, we propose to use basis functions of logarithmic growth,
%
\begin{equation}\label{eq:meixner log}
\psi_{nk}^{\mathrm{M},\log}(x) = P_k(\log x; t_n,\zeta).
\end{equation}
Then our convergence result for the Meixner process (Theorem \ref
{thm:main} below)
can be applied. Similarly, models based on the geometric Poisson (\cite
{Sh04}, Section 112.7.1)
or geometric Pascal processes can be reduced to
the linear case by modifying their respective basis functions analogously.

\section{Unbounded state space: The single-period problem}\label{se:single}

\subsection{Main result and first steps of the proof}

We now state our main result about the single-period problem where our
option has the exercise times $0=t_0<t_1<t_2$.
As noted above, the geometric Meixner model is contained in this result
by modifying
the basis functions according to (\ref{eq:meixner log}).
\begin{theorem}\label{thm:main}
Suppose $m=2$, that $S_t$ is a Meixner process and that the basis
functions are as
in the first line of Table \ref{ta:levy}. Put $(u,v)=(8,8)$.
If the number $N$ of paths and the number $K$ of basis functions satisfy
$N \geq K^{(u+\varepsilon)K}$ for some positive~$\varepsilon$, then
\[
\lim_{N\to\infty} \sup_{|\beta|=1} \expec[|\beta-\hat{\beta
}|^2] = 0.
\]
If $N \leq K^{(v-\varepsilon)K}$, then
\[
\lim_{N\to\infty} \sup_{|\beta|=1} \expec[|\beta-\hat{\beta
}|^2] =
\infty.
\]
For the standard Poisson, Gamma and Pascal processes, with their
respective basis
functions from Table \ref{ta:levy}, the same holds
if $(u,v)$ is replaced by $(10,4)$, $(8,8)$ and $(11,7)$, respectively.
\end{theorem}

The announced critical rate $\log N/\log\log N$ in Table \ref
{ta:comparison} then follows from
the fact that the solution of $N = K^{cK}$ satisfies $K\sim c^{-1} \log
N / \log\log N$ (see, e.g., de Bruijn \cite{deBr58}).

Looking at\vspace*{2pt} (\ref{eq:upper bound}) and (\ref{eq:lower
bound}) we begin the proof of Theorem \ref{thm:main} by bounding
$\|\Psi_1\|$ and $\|\Psi_1^{-1}\|$, defined by (\ref{eq:Psi}) and Table
\ref{ta:levy}. As in Section \ref{se:stentoft}, the letter $c$ denotes
various positive constants whose value is irrelevant.
\begin{lemma}\label{le:norm}
As $K\to\infty$, the values $\|\Psi_1\|$ and $\|\Psi_1^{-1}\|$ grow at
most exponentially in all four cases (Meixner, Poisson,
Gamma and Pascal), except for $\|\Psi_1^{\mathrm{P}}\|\leq c^K K^K$
and $\|\Psi_1^{\mathrm{Pa}}\|\leq c^K K^{2K}$.
\end{lemma}
\begin{pf}
The estimates for the Meixner, Poisson and Pascal cases are easy consequences
of the orthogonality relations (\ref{eq:orth charl})--(\ref{eq:orth
meix-poll})
and Stirling's formula.
It remains to deal with the Gamma case. The parameter $t-1$ in the
martingale property (\ref{eq:mart gamma})
is not quite compatible with the orthogonality relation (\ref{eq:orth
lag}) of the Laguerre polynomials.
But by the formula \cite{Sz75}
\[
L_k^{(\alpha-1)}(x) = L_k^{(\alpha)}(x) - L_{k-1}^{(\alpha)}(x)
\]
we obtain
%
\begin{equation}\label{eq: mix lag}
\expec[\psi_{1k}^{\mathrm{G}}(G_{t_1})\psi_{1l}^{\mathrm
{G}}(G_{t_1})] =
\cases{
-\pmatrix{k+t_1\cr k}, &\quad $k =l-1$, \vspace*{4pt}\cr
\dfrac{2k+t_1}{k+t_1}\pmatrix{k+t_1\cr k}, &\quad $k = l$, \vspace*{4pt}\cr
-\pmatrix{k+t_1-1\cr k-1}, &\quad $k = l+1$, \vspace*{4pt}\cr
0, &\quad $|k-l|\geq2$,}
\end{equation}
hence, $\Psi_1^{\mathrm{G}}$ is tridiagonal. Since (\ref{eq: mix lag})
grows only polynomially in $k$,
it is clear that so does $\|\Psi_1^{\mathrm{G}}\|$. As for the
inverse, note that $\Psi_1^{\mathrm{G}}$ is diagonally dominant
so that it suffices to bound the diagonal elements of $(\Psi
_1^{\mathrm
{G}})^{-1}$
(see Nabben \cite{Na99}, Theorem~3.1); 
note that the $\tau_k$ from that theorem are all equal to $1$ in our
situation.)
The diagonal elements $e_k$ of $(\Psi_1^{\mathrm{G}})^{-1}$ can be
computed recursively by \cite{Na99}
\[
e_{KK} = \frac{K}{K+t_1} \pmatrix{K+t_1-1\cr K-1}^{-1} \leq c^K
\]
and
\[
e_{k-1,k-1} = \frac{k+t_1}{k}\biggl( \frac{2k+t_1}{k+t_1} e_{k,k} -
e_{k+1,k+1}\biggr),\qquad 1\leq k<K.
\]
A straightforward backward induction shows that this implies
\[
|e_{kk}| \leq\bigl(4(t_1+1)\bigr)^{K-k+1} e_{KK}, \qquad 0\leq k<K,
\]
hence, $\|(\Psi_1^{\mathrm{G}})^{-1}\|$ grows at most exponentially
too.\vadjust{\goodbreak}
\end{pf}

We proceed to bound the fourth order moments appearing in (\ref
{eq:upper bound}). Using (\ref{eq:conn})
and the martingale relation (\ref{eq:mart}), we obtain
%
\begin{eqnarray} \label{eq:moments}
&&\expec[\psi_{2j}(S_{t_2})^2 \psi_{1k}(S_{t_1})^2] \nonumber\\
&&\qquad=
\expec\Biggl[\sum_{i=0}^{2j} d_{ji}(t_2) \psi_{2i}(S_{t_2}) \times\sum
_{s=0}^{2k} d_{ks}(t_1) \psi_{1s}(S_{t_1})\Biggr] \nonumber\\[-8pt]\\[-8pt]
&&\qquad= \sum_{i=0}^{2j} \sum_{s=0}^{2k} d_{ji}(t_2) d_{ks}(t_1) \expec[
\expec[\psi_{2i}(S_{t_2})\mid S_{t_1}] \psi_{1s}(S_{t_1})] \nonumber
\\
&&\qquad= \sum_{i=0}^{2j} \sum_{s=0}^{2k} d_{ji}(t_2) d_{ks}(t_1) \expec
[\psi
_{1i}(S_{t_1}) \psi_{1s}(S_{t_1})].\nonumber
\end{eqnarray}
The linearization coefficients $d_{ki}$ from the expansion (\ref{eq:conn})
are well-studied objects for various families of orthogonal
polynomials. They have combinatorial
interpretations in terms of (generalized) derangements, rook
polynomials and matching polynomials.
See Zeng \cite{Ze92} for on overview of these properties, explicit
formulas and many references.
Paraphrasing some of these formulas (\cite{Ze92}, Corollary 2) we have
%
\begin{eqnarray}\quad
\label{eq:d P}
d_{ki}^{\mathrm{P}}(t_n) &=& t_n^{2k-i} k!^2 i!
\sum_{s\geq0} \frac{t_n^s}{(s-k)!^2 (s-i)! (2k+i-2s)!},
\\
\label{eq:d G}
d_{ki}^{\mathrm{G}}(t_n) &=& 2^{2k+i} k!^2i!
\sum_{s\geq0} \frac{(t_n-1)_s }{4^s (s-k)!^2 (s-i)! (2k+i-2s)!},
\\
\label{eq:d Pa}
d_{ki}^{\mathrm{Pa}}(t_n) &=& (1+q)^{2k+i} k!^2i! \frac{(t_n)_k^2}{(t_n)_i}
\sum_{s\geq0} \frac{(t_n)_s (1+q)^{-2s}q^{-s}}{ (s-k)!^2 (s-i)!
(2k+i-2s)!}, \\
\label{eq:d M}
d_{ki}^{\mathrm{M}}(t_n) &=& (-2\cot\zeta)^{2k+i} k!^2i!
\sum_{s\geq0} \frac{(t_n)_s (1+(\cot\zeta)^{-2})^s}{4^s (s-k)!^2
(s-i)! (2k+i-2s)!}.
\end{eqnarray}
Here it is understood that $1/n!=0$ for $n$ a negative integer, as is
natural when extending the factorial by the Gamma function.
Therefore, the sums in (\ref{eq:d P})--(\ref{eq:d M}) run from
$s=\max\{
i,k\}$ to $s=k+\lfloor i/2 \rfloor$.

\subsection{Moment bounds in the Poisson case}
\begin{figure}

\includegraphics{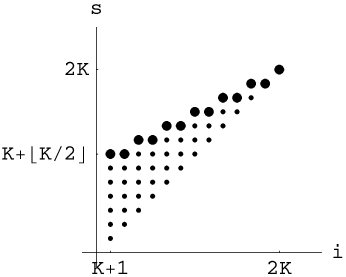}

\caption{The summation range of the first sum in
(\protect\ref{eq:est P 3}).}
\label{fig:sumrange}
\end{figure}

By (\ref{eq:orth charl}), (\ref{eq:moments}) and (\ref{eq:d P}), the
sum on the right-hand side of (\ref{eq:upper bound}) can be estimated
by
%
\begin{eqnarray}\label{eq:noname}
&&
\sum_{j=0}^K \sum_{k=0}^K \expec[\psi_{2j}^{\mathrm{P}}(S_{t_2})^2
\psi_{1k}^{\mathrm{P}}(S_{t_1})^2] \nonumber\\[-8pt]\\[-8pt]
&&\qquad\leq
c^K \sum_{j=0}^K \sum_{k=0}^K \sum_{i=0}^{2\min\{k,j\}} i!
\biggl( \sum_{s\geq0} b_{jis}^{\mathrm{P}} \biggr) \biggl( \sum
_{s\geq
0} b_{kis}^{\mathrm{P}} \biggr),\nonumber
\end{eqnarray}
where
\[
b_{kis}^{\mathrm{P}} := \frac{k!^2 i!}{(s-k)!^2 (s-i)! (2k+i-2s)!}.
\]
It is easy to see that $b_{k+1,i,k+l+1}^{\mathrm{P}} /
b_{k,i,k+l}^{\mathrm{P}} >1$ for $i\geq1$, $0\leq l\leq i/2$ and
$k\geq i-l$,
hence, $b_{k,i,k+l}^{\mathrm{P}}$ increases in $k$ under these
conditions. From this we deduce that the $s$-sums
in (\ref{eq:noname}) increase in $j$, respectively, $k$:
\begin{eqnarray*}
\sum_{s=\max\{i,k\}}^{k+\lfloor i/2 \rfloor} b_{kis}^{\mathrm{P}} &=&
\sum_{l=\max\{i-k,0\}}^{\lfloor i/2 \rfloor} b_{k,i,k+l}^{\mathrm
{P}} \\
&\leq&\sum_{l=\max\{i-k,0\}}^{\lfloor i/2 \rfloor}
b_{k+1,i,k+l+1}^{\mathrm{P}}\\
&=& \sum_{s=\max\{i,k\}+1}^{k+\lfloor i/2 \rfloor+1}
b_{k+1,i,s}^{\mathrm{P}}\\
&\leq&\sum_{s=\max\{i,k+1\}}^{k+\lfloor i/2 \rfloor+1}
b_{k+1,i,s}^{\mathrm{P}}.
\end{eqnarray*}
Using this in (\ref{eq:noname}) yields (recall that $c$ may change its
value in each occurrence)
%
\begin{eqnarray}\label{eq:est P}
&&
\sum_{j=0}^K \sum_{k=0}^K \expec[\psi_{2j}^{\mathrm{P}}(S_{t_2})^2
\psi
_{1k}^{\mathrm{P}}(S_{t_1})^2] \nonumber\\[-8pt]\\[-8pt]
&&\qquad\leq c^K K!^4 \sum_{i=0}^{2K} \Biggl( \sum_{s=K}^{K+\lfloor
i/2\rfloor
}\frac{i!^{3/2}}{(s-K)!^2(s-i)!(2K+i-2s)!} \Biggr)^2.
\nonumber
\end{eqnarray}
It is plain that the summand increases in $i$ for $K\geq0$, $0\leq
i\leq K$ and
$K\leq s\leq K+i/2$. Hence, we find that the portion $\sum_{i=0}^K$ of
the $i$-sum in (\ref{eq:est P}) can be bounded from above by
%
\begin{eqnarray}\label{eq:est P 2}
&&(K+1)K!^3 \Biggl( \sum_{s=K}^{\lfloor3K/2 \rfloor} \frac
{1}{(s-K)!^3(3K-2s)!} \Biggr)^2 \nonumber\\[-8pt]\\[-8pt]
&&\qquad\leq c^K K^{5K}.\nonumber
\end{eqnarray}
To see the\vspace*{1pt} last inequality, note that the summand in (\ref{eq:est P 2})
is unimodal with mode at
$s=K+K^{2/3}-\frac43 K^{1/3} + \Oh(1)$. Estimating this maximal
summand, by Stirling's formula
and some easy manipulations, shows that\vadjust{\goodbreak} the sum in (\ref{eq:est P 2})
is smaller than $c^K K^K$.
The remaining part $\sum_{i=K+1}^{2K}$ of the $i$-sum in (\ref{eq:est
P}) can be estimated by
%
\begin{eqnarray} \label{eq:est P 3}
&&\sum_{i=K+1}^{2K} i!  \Biggl( \sum_{s=i}^{K+\lfloor i/2\rfloor}
\frac
{i!s!}{(s-K)!^2(s-i)!(2K+i-2s)!} \Biggr)^2 \nonumber\\
&&\qquad\leq c^K \sum_{i=K+1}^{2K} i! \biggl( \frac{i!(K+\lfloor i/2 \rfloor)!}
{\lfloor i/2 \rfloor!^2(K+\lfloor i/2 \rfloor-i)!(i-2\lfloor i/2
\rfloor)!} \biggr)^2 \\
&&\qquad\leq c^K \sum_{i=K+1}^{2K} \frac{i!(K + \lfloor i/2 \rfloor)!^2}{(K
+\lfloor i/2 \rfloor-i)!^2} \nonumber
\leq c^K K^{6K}.
\end{eqnarray}
Note that in the first line we have introduced the new factor $s!$ in
the numerator. This makes the summand
increasing w.r.t. the substitution $i\to i+1$, $s\to s+1$. Hence, it
suffices to keep only the summands of the $s$-sum
with $s=K+\lfloor i/2 \rfloor$ (the thick dots in Figure \ref
{fig:sumrange}) which shows the first inequality.
%
As for the second inequality, note that
the factor $i!/\lfloor i/2 \rfloor!^2$ of the summand grows only
exponentially and that the factor
$(i - 2 \lfloor i/2 \rfloor)!$ in the denominator is clearly
negligible. Finally, the last sum in (\ref{eq:est P 3})
has increasing summands which, together with Stirling's formula,
implies the last inequality.
By (\ref{eq:est P}), the estimates (\ref{eq:est P 2}) and (\ref
{eq:est P 3}) show that
\[
\sum_{j=0}^K \sum_{k=0}^K \expec[\psi_{2j}^{\mathrm{P}}(S_{t_2})^2
\psi
_{1k}^{\mathrm{P}}(S_{t_1})^2] \leq c^K K^{10K}.
\]
In light of (\ref{eq:upper bound}) and Lemma \ref{le:norm}, the
value $u=10$ for
the Poisson process in Theorem~\ref{thm:main} is established.

As for the second assertion about the Poisson process in Theorem \ref
{thm:main}, note that, from (\ref{eq:moments}),
\[
\expec[\psi_{2K}(S_{t_2})^2 \psi_{1k}(S_{t_1})^2] =
\sum_{i=0}^{2K} \sum_{s=0}^{2k} d_{Ki}(t_2) d_{ks}(t_1) \expec[\psi
_{1i}(S_{t_1}) \psi_{1s}(S_{t_1})].
\]
The orthogonality property (\ref{eq:orth charl}) and formula (\ref
{eq:d P}) yield
\begin{eqnarray*}
\sum_{k=0}^K \expec[\psi_{2K}^{\mathrm{P}}(S_{t_2})^2 \psi
_{1k}^{\mathrm{P}}(S_{t_1})^2]
&\geq& c^K \sum_{k=0}^K \sum_{i=0}^{2k} d_{Ki}^{\mathrm{P}}(t_2)
d_{ki}^{\mathrm{P}}(t_1) i!\\[-1pt]
&\geq& c^K d_{K,2K}^{\mathrm{P}}(t_2) d_{K,2K}^{\mathrm{P}}(t_1) (2K)!
\\[-1pt]
&\geq& c^K (2K)!^3 \geq c^K K^{6K}.
\end{eqnarray*}
The second inequality follows from retaining only the summand $k=K$,
$i=2K$. This makes the sum
in (\ref{eq:d P}) collapse to the summand $s=2K$, hence, the third inequality.
Appealing to (\ref{eq:lower bound}) and Lemma \ref{le:norm} completes
the proof of the Poisson part
of Theorem \ref{thm:main}.
Note that the preceding estimates can presumably be improved.
This seems not worthwhile though; since our estimate of $\|\Psi
_1^{\mathrm{P}}\|$ in Lemma \ref{le:norm} is sharp,
we will not obtain equal values $u=v$ in Theorem \ref{thm:main} anyway,
unless at least one of the
bounds (\ref{eq:upper bound}) and (\ref{eq:lower bound}) was improved too.

\subsection{Moment bounds in the Meixner case}
The proofs in the remaining three cases are very similar to the Poisson
case. In the Meixner
case, we have
%
\begin{eqnarray}\label{eq:est M}
&&\sum_{j=0}^K \sum_{k=0}^K \expec[\psi_{2j}^{\mathrm{M}}(S_{t_2})^2
\psi
_{1k}^{\mathrm{M}}(S_{t_1})^2] \nonumber\\[-8pt]\\[-8pt]
&&\qquad\leq
c^K \sum_{j=0}^K \sum_{k=0}^K \sum_{i=0}^{2\min\{k,j\}} i!^2
\biggl( \sum_{s\geq0} b_{jis}^{\mathrm{M}} \biggr) \biggl( \sum
_{s\geq
0} b_{kis}^{\mathrm{M}} \biggr),\nonumber
\end{eqnarray}
where
\[
b_{kis}^{\mathrm{M}} := \frac{k!^2 s!}{(s-k)!^2 (s-i)! (2k+i-2s)!}.
\]
Again, $b_{k,i,k+l}^{\mathrm{M}}$ increases in $k$ and the remaining
steps to show the upper bound are
completely analogous to the Poisson case. This time the numerator
factor $s!$ in the analogue of (\ref{eq:est P 3})
appears naturally and is not introduced artificially to force some
monotonicity. Moreover, the lower
bound uses the same summands as in the Poisson case. Both resulting bounds
are of the form $c^K K^{8K}$, hence, $u=v=8$ in Theorem \ref{thm:main}.

\subsection{Moment bounds in the Pascal case}
We can reuse the values $b_{kis}^{\mathrm{M}}$ and the estimate that we
just sketched:
\begin{eqnarray*}
&&\sum_{j=0}^K \sum_{k=0}^K \expec[\psi_{2j}^{\mathrm
{Pa}}(S_{t_2})^2
\psi_{1k}^{\mathrm{Pa}} (S_{t_1})^2] \\[-1pt]
&&\qquad\leq c^K \sum_{j=0}^K \sum_{k=0}^K \sum_{i=0}^{2\min\{k,j\}} i!^3 k!
\biggl( \sum_{s\geq0} b_{jis}^{\mathrm{M}} \biggr) \biggl( \sum
_{s\geq
0} b_{kis}^{\mathrm{M}} \biggr) \\
&&\qquad\leq c^K K! (2K)! \sum_{j=0}^K \sum_{k=0}^K \sum_{i=0}^{2\min\{
k,j\}} i!^2
\biggl( \sum_{s\geq0} b_{jis}^{\mathrm{M}} \biggr) \biggl( \sum
_{s\geq
0} b_{kis}^{\mathrm{M}} \biggr) \\
&&\qquad\leq c^K K! (2K)! K^{8K} \leq c^K K^{11K}.
\end{eqnarray*}
The lower bound poses no new difficulties either.

\subsection{Moment bounds in the Gamma case}
This part is only slightly more involved. Due to (\ref{eq: mix lag}),
we have three $i$-sums
instead of one in the analogue of (\ref{eq:est P}). The right-hand side
of (\ref{eq: mix lag})
can be replaced by $c^K$ in each of these. Then one of the three
$i$-sums equals the $i$-sum in (\ref{eq:est M})
and the other two differ only in an index shift $b_{k,i\pm
1,s}^{\mathrm
{M}}$ which can be easily bounded by polynomial factors.
Thus the resulting growth rate is $c^K K^{8K}$, as for the Meixner
case. The proof of
Theorem \ref{thm:main} is complete.

\subsection{Side remark: The Bachelier model}
We finish this section with a remark about Brownian motion. If this is
the underlying process $S_t$, then appropriate basis functions can be
built from Hermite polynomials in such a way that $\|\Psi_1\|$,
$\|\Psi_1^{-1}\|$ and the analogue of (\ref{eq:orth charl}) grow only
exponentially \cite{GlYu04}. This is in line with the corresponding
growth orders in the Gamma and Meixner cases (and in the Poisson and
Pascal cases, if we renormalize our basis functions there by
$1/\sqrt{k!}$ and $1/k!$, resp.). What makes the Gaussian case peculiar
is that the linearization coefficients of the Hermite polynomials
induce only exponential growth too when plugged into~(\ref
{eq:moments}), whereas the linearization coefficients in the four cases
we treat in this paper grow faster.

\section{Unbounded state space: The multi-period problem}\label{se:multi}

In this section we extend the main result of the preceding section
(Theorem \ref{thm:main}) to the multi-period problem, that is, to $m+1$
exercise dates $0= t_0<\cdots< t_m$.
We know from the single-period problem that the critical rate cannot be
larger than
$\log N/\log\log N$, so we will be done if we can show that there is
an upper bound for the mean square error
of the form $K^{cK}$ for some positive $c$.
Fortunately, this can be deduced with little effort from a result of
Glasserman and Yu \cite{GlYu04}
and the estimates from the preceding section about the single-period problem.
Following \cite{GlYu04}, we assume that a representation analogous
to (\ref{eq:repr}) holds at time $t_m$ and
that the payoff functions do not grow too fast in the following sense.
\begin{theorem}\label{thm:mult}
Suppose that the payoff functions satisfy the growth constraint
\[
\expec[h_n(S_{t_n})^4] \leq\max_\nu\biggl(\frac{t_{\nu+1}}{t_\nu
}\biggr)^{2K}
\max_{\nu,k}\expec[\psi_{\nu k}(S_{t_\nu})^4],\qquad 0\leq n\leq m.
\]
Then the mean square error of the estimated coefficients satisfies
\[
\sup_{|\beta_{m-1}|=1} \expec[|\beta_n-\hat{\beta}_n|^2] \leq N^{-1}
c^K K^{(m-n+1)uK}, \qquad 1\leq n<m,
\]
where $u$ takes on the same values as in Theorem \ref{thm:main},
that is, $8, 10, 8, 11$ for $S_t$,
the Meixner, standard Poisson, Gamma and Pascal process, respectively.
\end{theorem}
\begin{pf}
By results of Glasserman and Yu \cite{GlYu04}, Theorem 3 and the last
formula before (18) on page 2096 and Jensen's inequality, we have
\begin{eqnarray*}
&&\sup_{|\beta_{m-1}|=1} \expec[|\beta_n  - \hat{\beta}_n|^2] \\
&&\qquad\leq\frac{c^K}{N} \max_{1\leq\nu<m} \| \Psi_{\nu}^{-1}\|^3 \max
_{\nu,k} \expec[\psi_{\nu k}(S_{t_\nu})^4]^{m-n}
\max_{\nu,k} \expec[\psi_{\nu k}(S_{t_\nu})^2]^2 \\
&&\qquad\leq\frac{c^K}{N} \max_{1\leq\nu<m} \| \Psi_{\nu}^{-1}\|^3 \max
_{\nu
,k} \expec[\psi_{\nu k}(S_{t_\nu})^4]^{m-n+1}.
\end{eqnarray*}
Note that Glasserman and Yu \cite{GlYu04} assume that the
moments $\expec[\psi_{n k}(S_{t_\nu})^2]$ and $\expec[\psi_{n
k}(S_{t_\nu})^4]$
are increasing in $n$ and $k$ and formulate their Theorem 3
with $\expec
[\psi_{mK}^{2(4)}]$ instead of
$\max_{\nu,k}\expec[\psi_{\nu k}^{2(4)}]$. But an inspection of their
proof quickly shows that taking the max
in the above estimate gets rid of the monotonicity assumption.
Now note that $ \| \Psi_{\nu}^{-1}\|\leq c^K$ in all our four cases by
Lemma \ref{le:norm} and that
\[
\max_{\nu,k} \expec[\psi_{\nu k}(S_{t_\nu})^4] \leq\max_\nu
\sum_{j=0}^K \sum_{k=0}^K \expec[\psi_{\nu j}(S_{t_\nu})^2 \psi
_{\nu
k}(S_{t_\nu})^2] \leq c^K K^{uK},
\]
where the double sum has been estimated in the proof of Theorem \ref{thm:main}.
\end{pf}

We have thus seen that Table \ref{ta:comparison} correctly describes
the general (i.e., multi-period) situation.

\section{Conclusion}
Stentoft \cite{St04} and Glasserman and Yu \cite{GlYu04} obtained
apparently contradictory results about the convergence of the
Longstaff--Schwartz algorithm. The main difference between their
respective assumptions is
the (un-)boundedness of the support of the underlying at the exercise
dates. In this light
the pessimistic results of Glasserman and Yu (and our Theorems \ref
{thm:main} and \ref{thm:mult}) turn out to stem from the tails
of the distribution of the underlying.

The present paper shows that Stentoft's result can be applied to L\'evy
models under mild assumptions
and extends Glasserman and Yu's \cite{GlYu04} results to several
concrete processes.
Thus we provide some evidence that Glasserman and Yu \cite{GlYu04} were
right to
conjecture that their results for Brownian motion and geometric
Brownian motion
extend to other models.

Concerning Stentoft \cite{St04} and our Section \ref{se:stentoft}:
although the boundedness of the underlying induces a nice (polynomial)
relation between the
number of basis functions and the necessary number of Monte Carlo paths,
it seems not yet completely clear that it is a harmless assumption in practice.
A natural question for future research is how strongly the size of the
truncation intervals influences
the convergence speed of the calculated prices.

\section*{Acknowledgments}
I thank Lars Stentoft, Friedrich Hubalek and an anonymous referee
for helpful comments.

%

%
\printaddresses

\end{document}